\documentclass[letter,11pt,leqno, english]{amsart}
        \title[Rational magnetic equivariant K-theory]
				{Rational magnetic equivariant K-theory}
        \author{Higinio Serrano}
      	\address{Departamento de Matem\'aticas, Centro de Investigaciones y Estudios Avanzados, Av. Instituto Politécnico Nacional 2508,
Col. San Pedro Zacatenco, Mexico D.F.,CP 07360. Mexico}
	       \email{hserrano@math.cinvestav.mx}
 
     \author{Bernardo Uribe}
      \address{Departamento de Matem\'aticas y Estad\'istica\\
       Universidad del Norte\\ 
       Km. 5 via Puerto Colombia, Barranquilla, Colombia}
      \email{bjongbloed@uninorte.edu.co, buribe@gmail.com}
 \author{Miguel A. Xicotencatl}
      	\address{Departamento de Matem\'aticas, Centro de Investigaciones y Estudios Avanzados, Av. Instituto Politécnico Nacional 2508,
Col. San Pedro Zacatenco, Mexico D.F.,CP 07360. Mexico}
	       \email{xico@math.cinvestav.mx}
            \keywords{Magnetic group; Equivariant K-theory; Real K-theory; Queternionic K-theory}       
       \subjclass[2020]{19L47,19L50,55N91}
\thanks{}
\date{\today}
\usepackage{amssymb}
\usepackage{pdfsync}
\usepackage{hyperref}
\usepackage{calc}
\usepackage{enumerate,amssymb,amscd}
\usepackage{color}
\usepackage[arrow,curve,matrix,tips,2cell]{xy}
  \SelectTips{eu}{10} \UseTips
  \UseAllTwocells
\usepackage{tikz}
\usepackage{pdfcolmk}

\DeclareMathAlphabet{\matheurm}{U}{eur}{m}{n}

\makeatletter
\newcommand{\colim@}[2]{%
  \vtop{\m@th\ialign{##\cr
    \hfil$#1\operator@font colim$\hfil\cr
    \noalign{\nointerlineskip\kern1.5\ex@}#2\cr
    \noalign{\nointerlineskip\kern-\ex@}\cr}}%
}
\newcommand{\colim}{  \mathop{\mathpalette\colim@{\rightarrowfill@\textstyle}}\nmlimits@}
\makeatother

\newcounter{commentcounter}

\theoremstyle{plain}

\newtheorem*{theorem*}{Theorem}
\newtheorem*{mtheorem*}{Main Theorem}

\theoremstyle{definition}

\theoremstyle{remark}
\newtheorem*{summary*}{Summary}
\newtheorem{innercustomgeneric}{\customgenericname}
\providecommand{\customgenericname}{}
\newcommand{\newcustomtheorem}[2]{%
  \newenvironment{#1}[1]
  {%
   \renewcommand\customgenericname{#2}%
   \renewcommand\theinnercustomgeneric{##1}%
   \innercustomgeneric
  }
  {\endinnercustomgeneric}
}

\newcustomtheorem{customthm}{Theorem}
\newcustomtheorem{customcor}{Corollary}

\makeatletter\let\c@equation=\c@theorem\makeatother

\newcommand{\version}[1] 
{\begin{center} Last edited on #1\\
    Last compiled on \today\\
file name: \jobname
  \end{center}
}

\begin{document}

\begin{abstract}  
We introduce the magnetic equivariant K-theory groups as the K-theory groups associated to magnetic groups
and their respective magnetic equivariant complex
bundles. We restrict the magnetic group to its subgroup of elements that act complex linearly, and we show that
this restriction induces a rational isomorphism with the conjugation invariant part of the complex
equivariant K-theory of the restricted group. This isomorphism allows to calculate the torsion free
part of the magnetic equivariant K-theory groups reducing it to known calculations in complex
equivariant K-theory
 \end{abstract}

\maketitle

\section*{Introduction}

The discovery \cite{Thouless}, prediction \cite{Haldane}  and later experimental observation \cite{Observation_QAHE} of the Quantum Anomalous Hall effect in magnetic topological insulators has been a recent and very exciting development in the realm of
condensed matter physics.
One of the key features of this phenomenon is the fact that it is quantized; namely the Hall conductivity can only acquire values proportional to integer multiples of the von Klitzing constant $(e^2/ \hbar)$. The integer here is the Chern number of the
vector bundle of valence states, or Bloch bundle, of the material. More remarkably, this Chern number was also shown
to determine a quantized Hall conductivity on materials that need
to be modeled with tools of noncommutative geometry \cite{Bellissard1, Bellissard2}.

Later, the proposal for the existence of a Quantum spin Hall Effect, namely a quantized spin-Hall conductance with vanishing charge-Hall conductance, was put forward \cite{KaneMele} \cite{Bernevig} and was experimentally observed \cite{Observation_QSHE}. Here the invariant which characterizes a state as trivial or non-trivial band insulator is a number in the group $\mathbb{Z}_2$ (regardless if the state exhibits or does not exhibit a quantum spin Hall effect). This number is not zero
whenever the valence bands generate the $\mathbb{Z}_2$-invariant
which was shown to live in Atiyah's real K-theory of the 2D torus \cite{KaneMeleZ2}. This invariant is known as the Kane-Mele invariant.

Both phenomena are shown to be related to the topological
properties of the Bloch bundle of valence states of the material.
The first Chern number being the integer invariant in magnetic
2D materials, and the Kane-Mele invariant being the torsion invariant in materials that preserve time reversal symmetry.
Among the many interesting properties these quantized quantities infer in the electronic properties of a material is that
they are robust under adiabatic perturbations of the Hamiltonian; in other words, the effects are robust under the presence of small impurities in the materials \cite{QAHE-impurities}. 

These new electronic phases in materials are nowadays the subject of intense research and have opened several new roads for the discovery and classification of new compounds. The importance
and relevance of these new ``Topological Phases of Matter" was
recognized in the Nobel Prize in Physics of 2016 which was awarded to David J. Thouless,  F. Duncan M. Haldane and J. Michael Kosterlitz ``for their theoretical discoveries of topological phase transitions and topological phases of matter" \cite{Nobel}.

The understanding of the topological invariants of crystals, magnetic or not, depends on the explicit knowledge of the
group of symmetries of the crystal, as well as 
the equivariant K-theory groups of the 2D and 3D torus.
The group of symmetries that are of interest incorporate the
group of spatial symmetries of the crystal, as well as the 
symmetries which composed with the time reversal operator preserves
the Hamiltonian. These groups have been known in the physics literature as magnetic groups or \v Subnikov groups \cite{Subnikov}, while in the mathematics community they have
been simply known as $\mathbb{Z}_2$-groups. In this paper we take
the physical name of magnetic groups and we define
their associated equivariant K-theory groups. 
To differentiate them from the well-known complex equivariant K-theory groups \cite{Atiyah-Segal} we have coined them {\it Magnetic Equivariant K-theory groups}.

The magnetic equivariant K-theory groups have been studied as early as the year 1970 \cite{Karoubi}, have been recently articulated to the
to the study of topological phases of matter \cite{Freed-Moore}, \cite{Gomi}, and have been further developed to provide tools for their explicit calculation \cite{Gomi-AHSE}, \cite{Gomi-TCM}, \cite{Serrano}.

The magnetic equivariant K-theory groups are in general not easy to calculate. Each case needs to be treated separately, and an appropriate cell subdivision of the underlying space is necessary for applying the cohomological tools that help determine the desired groups. The K-theory groups may be torsion and non-torsion, but sometimes the relevant feature appears as a non-torsion invariant. In those particular cases the rational magnetic equivariant K-theory is enough to understand the non-torsion invariants. This is the subject of this work.

We first take the restriction map from the magnetic equivariant K-theory
to the complex equivariant K-theory of the underlying subgroup that acts complex linearly. We show that its image lies in the conjugation invariant subgroup of this equivariant complex
 K-theory, and moreover, that this restriction map induces a rational isomorphism. This is our key result and it is the main theorem of the paper.

We further show that the rational isomorphism also applies to the twisted version of the magnetic equivariant K-theory, and we finish with an application to 2D materials which preserve
the combination of a four-fold symmetry and time reversal, as well as the spin in the $z$ direction. Here we show that the
Chern number of the spin up valence bands is indeed an integer invariant, and using results from \cite{C4T-Hamiltonian}, we show that the parity of this Chern number provides the $\mathbb{Z}_2$-invariant in 2D topological insulator altermagnets.

\section{Magnetic equivariant K-theory}

Magnetic equivariant K-theory is the K-theory of complex vector bundles
with actions of magnetic groups \cite{Karoubi},\cite{Freed-Moore}, \cite{Serrano}. Let us be more precise.

A {\it{magnetic group}}  consist of a group $G$ together with a surjective homomorphism to $\mathbb{Z}_2$. Denote this map $\phi: G \to \mathbb{Z}_2$ and call $G_0$ the kernel of $\phi$, thus obtaining the following short exact sequence:
\begin{align}
1 \to G_0 \to G \stackrel{\phi}{\to} \mathbb{Z}_2 \to 1.
\end{align}
A homomorphism of magnetic groups $f: G \to L$ is a homomorphism of groups
which is compatible with the surjective maps to $\mathbb{Z}_2$. A subgroup $H \subset G$ of a magnetic group is also a magnetic group if $H \cap G_0 \neq H$;
otherwise $H \subset G_0$ and the subgroup is not magnetic (just a group). 

Examples of magnetic groups are among others, magnetic space groups in crystallography \cite{Heesch} and \v Subnikov groups \cite{Subnikov}.

Here we will restrict to the case on which the magnetic group $G$ is also a compact Lie group. For the applications to crystallography $G$ will be a magnetic point group.

Let $X$ be a compact $G$ space and define a {\it magnetic equivariant vector bundle}
as a complex vector bundle $E \stackrel{p}{\to} X$ endowed with an action of $G$ compatible with the map $p$ such that $G_0$ acts complex linearly on the fibers and $G \backslash G_0$ acts complex anti-linearly.

A homomorphism of magnetic equivariant vector bundles is simply a $G$-equivariant
homomorphism of the underlying complex vector bundles. We may take the isomorphism classes of magnetic equivariant vector bundles and we may define the {\it magnetic $G$-equivariant K-theory of $X$} as the Grothendieck group of the isomorphism classes. We will denote this K-theory with a calligraphic letter $\mathcal{K}_G(X)$, consisting of formal differences $E_0-E_1$ of magnetic equivariant bundles subject to the equivalence relation
\begin{align}
    E_0-E_1 \sim E_0'-E_1' \Longleftrightarrow E_0 \oplus E_1' \oplus F \cong E_0' \oplus E_1 \oplus F
\end{align}
for some  magnetic $G$-equivariant vector bundle $F$. For non-magnetic groups,
such as $G_0$, we will denote  $K_{G_0}(X)$ the complex $G_0$-equivariant
K-theory groups of $X$ \cite{Segal}.

A $G$-equivariant map $\psi : X \to Y$ induces a homomorphism of abelian groups 
$\psi^* : \mathcal{K}_G(Y) \to \mathcal{K}_G(X)$, $E \mapsto \psi^*E$, making
 $\mathcal{K}_G$ a contravariant functor from compact $G$-spaces to abelian groups. A homomorphism of magnetic groups $\alpha : H \to G$ induces
a homomorphism $\mathcal{K}_G(X) \to \mathcal{K}_H(X)$, and  when $H$
is non-magnetic and $\alpha$ factors through $H \to G_0 \subset G$ the homomorphism is $\mathcal{K}_G(X) \to {K}_H(X)$. 

For $G$-spaces $X$ with a choice of base point $x_0 $, the reduced K-theory groups are defined as the kernel of the pullback under the restriction map:
\begin{align}
\widetilde{\mathcal{K}}_G(X) : = \ker \left( \mathcal{K}_G(X) \to \mathcal{K}_G(\{x_0\} 
\right).
\end{align}
For a pair $(X,Y)$ of compact $G$-spaces we define $\mathcal{K}_G(X,Y) : = \widetilde{\mathcal{K}}_G(X/Y)$ and whenever $Y = \emptyset$ we set $X/Y:=X_+$,
thus having the isomorphism $\mathcal{K}_G(X,\emptyset) \cong \mathcal{K}_G(X)$.

The higher K-theory groups for $q \in \mathbb{N}$ are defined as follows:
\begin{align}
    \mathcal{K}^{-q}_G(X,Y) := \mathcal{K}_G(X \times B^q, X \times S^{q-1}\cup Y \times B^q), \label{suspension}
\end{align}
thus having the usual suspension equality $\widetilde{\mathcal{K}}_G^{-q}(X) =
\widetilde{\mathcal{K}}_G(\Sigma^q X)$.

The magnetic $G$-equivariant K-theory groups are one example of a $G$-equivariant cohomology \cite{May}. It was firstly defined in \cite{Karoubi}, and further elaborated in \cite{Freed-Moore} and \cite{Gomi}. The present description has been developed by the first author in \cite{Serrano} and its structural properties and applications will be presented in a forthcoming publication.

\subsection*{Properties}
Among the many properties that the magnetic $G$-equivariant K-theory groups have, we want to highlight the following:
\begin{itemize}
    \item The magnetic $G$-equivariant K-theory over a point $x_0$ is the group 
    of isomorphism classes of magnetic representations of $G$. These representations were called {\it corepresentations} by Wigner and many of their properties were described on his book \cite{Wigner}.
    A magnetic representation of $G$ is a complex vector space $V$ with a $G$ action that is complex linear on $G_0$ and complex anti linear on $G \backslash G_0$. The abelian group of isomorphism classes of magnetic representation is a free $\mathbb{Z}$-module generated by the irreducible ones. Denote this group with calligraphic letter $\mathcal{R}(G)$ and we have
    \begin{align}
        \mathcal{K}_G(\{x_0\}) \cong \mathcal{R}(G).
    \end{align}

    The restriction of the magnetic representations to $G_0$ defines a homomorphism 
    \begin{align}
        \mathcal{R}(G) \to R(G_0), \ \ V \mapsto V|_{G_0}
    \end{align}
    to the abelian group $R(G_0)$ of isomorphism classes of complex $G_0$-representations. Any irreducible magnetic representation $V$ of $G$  decomposes into irreducible representations of $G_0$  
    in $R(G_0)$ fitting only one of the following three cases:
   \begin{itemize}
       \item Real type: $V|_{G_0} \cong U$ with $U$ irreducible $G_0$-representation.
       \item Complex type: $V|_{G_0} \cong W \oplus \widehat{W}$ with $W$ irreducible, $\widehat{W}$ the conjugate representation defined in eqn. \eqref{conjugate representation}, and $W \ncong \widehat{W}$ as $G_0$-representations.
        \item Quaternionic type: $V|_{G_0} \cong Z \oplus Z$ with $Z$ irreducible $G_0$-representation. Here we have that $Z \cong \widehat{Z}$.
   \end{itemize}
   Hence, the abelian group $\mathcal{R}(G)$ may be split as 
   \begin{align}
       \mathcal{R}(G) \cong \mathcal{R}(G)_{\mathbb{R}}\oplus \mathcal{R}(G)_{\mathbb{C}}\oplus \mathcal{R}(G)_{\mathbb{H}}
   \end{align}
   where $\mathcal{R}(G)_{\mathbb{F}}$ corresponds to the magnetic representations of type $\mathbb{F}$ with $\mathbb{F}$ one of the commuting fields $\mathbb{R}$,$\mathbb{C}$, or the division ring $\mathbb{H}$, respectively.
\item Whenever $X$ is a trivial $G$-space, there is a canonical decomposition
\begin{align}
    \mathcal{K}_G(X) \cong \left( \mathcal{R}(G)_{\mathbb{R}} \otimes KO(X) \right) \oplus  \left(\mathcal{R}(G)_{\mathbb{C}}\otimes K(X)  \right)  \oplus  \left( \mathcal{R}(G)_{\mathbb{H}}\otimes KSp(X) \right)
\end{align}
   where $KO(X)$, $K(X)$ and $KSp(X)$ denote the Grothendieck groups of real, complex and quaternionic vector bundles over $X$. If $V$ is an irreducible magnetic representation of $\mathbb{F}$-type and $E$ is a magnetic $G$-equivariant vector bundle, the bundle $\mathrm{Hom}(V,E)$ defines an element in 
   K-theory of $\mathbb{F}$-type. Carrying out this assignment for all irreducible magnetic representations the isomorphism of above follows.
\item The coefficients of the magnetic $G$-equivariant cohomology theory split as follows. 
If * denotes the one-point trivial $G$-space, there is a canonical decomposition:
\begin{align}
    \mathcal{K}^{-q}_G(*) \cong \left( \mathcal{R}(G)_{\mathbb{R}} \otimes KO^{-q}(*) \right) \oplus  \left(\mathcal{R}(G)_{\mathbb{C}}\otimes K^{-q}(*)  \right)  \oplus  \left( \mathcal{R}(G)_{\mathbb{H}}\otimes KSp^{-q}(*) \right), \label{splitting suspension}
\end{align}
    where $KO^*$, $K^*$ and $KSp^*$ denote the K-theory of real, complex and quaternionic vector bundles. 
   This result is proven by the first author in \cite{Serrano} and generalizes
   the statement once restricted to real $G$-equivariant K-theory \cite{Atiyah-Segal}.
\item The magnetic $G$-equivariant K-theory is 8 periodic \cite{Gomi, Serrano}:
\begin{align}
    \mathcal{K}_G^{-q-8}(X) \cong \mathcal{K}_G^{-q}(X).
\end{align}
    One can therefore define the positively graded magnetic equivariant K-theory groups using this periodicity.  
\item The restriction to orbit types gives the following isomorphisms. For $H \subset G$ we have:
\begin{align}
    \mathcal{K}^*_G(G/H) \cong \left\{
    \begin{matrix}
    \mathcal{K}^*_H(*) &  \mathrm{if} & H \not \subset G_0 \\
    K_H^*(*) &  \mathrm{if} & H \subset G_0. 
    \end{matrix}
    \right.
\end{align}
    Moreover, if $N \subset G_0$ is a normal subgroup $N \leq G$, and $N$ acts
    freely on $X$, then the projection $\pi : X \to X/N$ induces an isomorphism
    \begin{align}
        \pi^*:\mathcal{K}_{G/N}^*(X/N) \stackrel{\cong}{\to} \mathcal{K}_{G}^*(X). 
    \end{align}
\end{itemize}

\subsection*{Calculation}
Given a magnetic group $G$ and a $G$-space, the calculation of the magnetic $G$-equivariant K-theory groups is not straightforward. Perhaps the most common
form to calculate these K-theory groups is using a $G$-CW decomposition and the spectral sequence that the decomposition induces. The first two pages of this
spectral sequence are manageable, but the extension problems that the higher
differentials encode makes this procedure difficult and hard to use for non experts.

The magnetic $G$-equivariant K-theory groups possess both torsion and non-torsion information. Sometimes most of the non-trivial information is torsion, but in some others, the non-torsion part is already good enough.

In what follows we will outline a procedure for extracting the non-torsion
information of the magnetic $G$-equivariant K-theory groups. The idea is 
to study the restriction map to the complex $G_0$-equivariant K-theory and to 
determine its image.

\section{Rational magnetic equivariant K-theory}

Denote by $\iota: G_0 \to G$ the inclusion homomorphism of $G_0 = \mathrm{Ker}(\phi)$ 
into $G$. The restriction homomorphism 
\begin{align}
    \iota^*:\mathcal{K}_G(X) \to K_{G_0}(X)
\end{align}
maps magnetic $G$-equivariant vector bundles to complex $G_0$-equivariant ones \cite{Segal}.

We claim that there is an action  of $\mathbb{Z}_2$ on $K_{G_0}(X)$
such that the image of $\iota^*$ lands in the $\mathbb{Z}_2$ invariant part. Let us see how this action is defined.

Take $F \stackrel{p}{\to} X$ a complex $G_0$-equivariant bundle over $X$. Choose any element
$a_0 \in G \backslash G_0$ and define the pullback bundle
\begin{align}
    a_0^{*}\overline{F} = \{(\overline{s},x) \in \overline{F}\times X | p(s)=a_0x \}
\end{align}
where $\overline{F}$ denotes the complex conjugate bundle of $F$. Endow
$a_0^{*}\overline{F}$ with the $\mathbb{C}$-module structure as follows: for
$\lambda \in \mathbb{C}$ let $\lambda \cdot (\overline{s},x) := (\overline{\lambda s},x)$, and endow $a_0^{*}\overline{F}$ with the following
$G_0$-equivariant structure: for $g \in G_0$, $g \cdot (\overline{s},x) = (\overline{a_0ga_0^{-1}s},gx)$. Note that $(\overline{a_0ga_0^{-1}s},gx)$ belongs
to $a_0^{*}\overline{F}$ since $p(a_0ga_0^{-1}s)=a_0ga_0^{-1}p(s)=a_0gx$.

We have now that $a_0^{*}\overline{F}$ is a complex $G_0$-equivariant bundle
over $X$. Applying the construction again, we get 
\begin{align}
a_0^*\overline{(a_0^{*}\overline{F})} = (a_0^2)^*F.
\end{align}
Since  $a_0^2$ belongs to $G_0$, we have that $F$ and $(a_0^2)^*F$ become
isomorphic $G_0$-equivariant bundles with the following homomorphism:
\begin{align}
    F \stackrel{\cong}{\to} (a_0^2)^*F, \ \ s \mapsto (a_0^2 s, p(s)). 
\end{align}

We have therefore defined an involution
\begin{align}
  K_{G_0}(X) \to K_{G_0}(X),  \ \ \ F \mapsto a_0^*\overline{F}
\end{align}
which makes $K_{G_0}(X)$ a $\mathbb{Z}_2$-module. 
Note that this involution is independent of the choice of element in $G \backslash G_0$ since $a_0^{*}\overline{F}$ and $(ga_0)^{*}\overline{F}$
are isomorphic for any $g \in G_0$. The homomorphism
\begin{align}
    a_0^{*}\overline{F} \stackrel{\cong}{\to} (ga_0)^{*}\overline{F}, \ \ \ (\overline{s},x) \mapsto (\overline{ga_0s},x)
\end{align}
gives the desired $G_0$-equivariant isomorphism (here $s \in F$ and $p(s)=a_0x$).

In the case that $X$ is a point, the $\mathbb{Z}_2$ action on $R(G_0)$ provides the conjugate representation that was mentioned above. The involution is then:

\begin{align}
  R(G_0) \to R(G_0), \ \ \ W\mapsto \widehat{W}:=a_0^* \overline{W}.
\label{conjugate representation}
\end{align}

Now we are ready to state the main result of this work.

{\bf Theorem.} {\it Let $X$ be a compact $G$-space, $G \stackrel{\phi}{\to} {\mathbb{Z}_2}$ a magnetic group and
$\iota : G_0 \to G$ the inclusion of the kernel of $\phi$.
Then the pullback of the restriction $\iota^*: \mathcal{K}_G^*(X) \to K_{G_0}^*(X)$ lands in the $\mathbb{Z}_2$-invariant subgroup}
\begin{align}
    \iota^*: \mathcal{K}_G^*(X) \to K_{G_0}^*(X)^{\mathbb{Z}_2}
\end{align}
{\it and it becomes an isomorphism rationally}
\begin{align}
    \iota^* \otimes \mathbb{Q}: \mathcal{K}_G^*(X) \otimes \mathbb{Q} \stackrel{\cong}{\to} K_{G_0}^*(X)^{\mathbb{Z}_2}\otimes \mathbb{Q}.
\end{align}

{\bf Proof.}
First let us show that the image of $\iota^*$ lands in the 
$\mathbb{Z}_2$ invariant part. Take a magnetic $G$-equivariant
vector bundle $E \stackrel{p}{\to} X$ and denote its restriction $\iota^*E$ to $G_0$ with the same letter $E$. The $\mathbb{Z}_2$
action takes $E$ to its conjugate $a_0^*\overline{E}$
where $a_0 \in G \backslash G_0$. The homomorphism
\begin{align}
    E \to a_0^*\overline{E},  \ \ \ s \mapsto (\overline{a_0s},p(s))
\end{align}
is the desired $G_0$-equivariant isomorphism. Hence
$\iota^*E \in K_{G_0}^*(X)^{\mathbb{Z}_2}$.

Now, since the homomorphism $\iota^*$ is a natural transformation of functors, it is enough to show that it induces
an isomorphism at the level of orbit types. Let us show this.

Take $H\stackrel{\varphi}{\to}{\mathbb{Z}_2} $ any magnetic group with $H_0= \mathrm{ker}(\varphi)$ and note that the magnetic irreducible representations of $H$, as well as the complex irreducible representations of $H_0$, may be broken
into real, complex and quaternionic representations \cite{Wigner}; the former one depending on the automorphism
that $H/H_0$ induces on $H_0$. Split the magnetic representations of $H$ and the complex representations of $H_0$ accordingly:
\begin{align}
\mathcal{R}(H) =& \mathcal{R}(H)_{\mathbb{R}}  \oplus
\mathcal{R}(H)_{\mathbb{C}} \oplus \mathcal{R}(H)_{\mathbb{H}}\\ 
{R}(H_0) =& {R}(H_0)_{\mathbb{R}}  \oplus
{R}(H_0)_{\mathbb{C}} \oplus {R}(H_0)_{\mathbb{H}}
\end{align}

A complex irreducible representation $V$ of $H_0$ is of complex type with respect to $H$ if the conjugate representation $\widehat{V}:=a_0^* \overline{V}$ is not isomorphic to $V$ (here $a_0$ is any element in $H \backslash H_0$). Notice that in this case $\widehat{V}$ is simply $\overline{V}$, and therefore we will use this identification $\widehat{V}:=\overline{V}$ hereafter. Whenever $\overline{V} \cong V$, let $T \in \mathrm{Hom}_{\mathrm{Rep}(H_0)}( \overline{V}, V)$ be
the isomorphism and denote $\rho:H_0 \to \mathrm{GL}(V)$ and
$\rho':H_0 \to \mathrm{GL}( \overline{V})$ (here  $\rho'(h)=\overline{\rho(a_0ha_0^{-1})}$) the homomorphisms corresponding to the complex representations $V$ and $\overline{V}$.
We have then
\begin{align}
    \rho'=T^{-1} \rho T  \  \  \ \mathrm{and} \  \  \ \rho(a_0^{-2}) = \pm T \overline{T}.
\end{align}
The complex representation $V$ is called of real type if  
$\rho(a_0^{-2}) = + T \overline{T}$ and of quaternionic type if $\rho(a_0^{-2}) = - T \overline{T}$. 

The restriction homomorphism $\mathcal{R}(H) \to R(H_0)$
splits into three maps $\mathcal{R}(H)_{\mathbb{F}} \to R(H_0)_{\mathbb{F}}$ for $\mathbb{F}$ in $\{\mathbb{R}, \mathbb{C},\mathbb{H} \}$. In the real case the homomorphism
\begin{align}
    \mathcal{R}(H)_{\mathbb{R}} \stackrel{\cong}{\to} R(H_0)_{\mathbb{R}}
\end{align}    
is an isomorphism with $R(H_0)_{\mathbb{R}}$ fixed by the $\mathbb{Z}_2$ action. In the quaternionic case the homomorphism
\begin{align}
    \mathcal{R}(H)_{\mathbb{H}} \to R(H_0)_{\mathbb{H}} \  \ \ \ \ \  \  \mathcal{R}(H)_{\mathbb{H}}  \otimes \mathbb{Q} \stackrel{\cong}{\to} R(H_0)_{\mathbb{H}} \otimes \mathbb{Q}
\end{align}   
is injective and of full rank (it is multiplication by 2 on the generators), and therefore an isomorphism rationally, with  $R(H_0)_{\mathbb{H}}$ fixed by $\mathbb{Z}_2$.
In the complex case the homomorphism
\begin{align}
    \mathcal{R}(H)_{\mathbb{C}} \stackrel{\cong}{\to} R(H_0)_{\mathbb{C}}^{\mathbb{Z}_2}
\end{align} 
is an isomorphism with the $\mathbb{Z}_2$ invariant part.

Hence the restriction homomorphism
\begin{align}
    \mathcal{R}(H) \to R(H_0)^{\mathbb{Z}_2} 
\end{align} is injective and of full rank. Rationally we get an isomorphism
\begin{align}
    \mathcal{R}(H) \otimes \mathbb{Q} \to R(H_0)^{\mathbb{Z}_2} \otimes \mathbb{Q}. 
\end{align}

Now, let us take any $G$-orbit type $G/H$ for $H \subset G$. We have two cases, either $H \subset G_0$ or not. Assuming that
$H \subset G_0$ we have that
\begin{align}
\mathcal{K}_G(G/H) \cong K_H(*) \ \ \ \mathrm{and} \ \ 
{K}_{G_0}(G/H)\cong K_{H\times G_0}(G) \cong K_{H}(*) \oplus K_{H}(*)
\end{align} 
where the $\mathbb{Z}_2$ action on $K_{H}(*) \oplus K_{H}(*)$ 
swaps the summands. We have that the restriction map in this case induces an isomorphism
\begin{align}
\xymatrix{
    \mathcal{K}_G(G/H) \ar[r] \ar[d]^\cong & K_{G_0}(G/H)^{\mathbb{Z}_2} \ar[d]^\cong \\
    K_{H}(*) \ar[r]^{\cong \hspace{1cm}} & \left(K_H(*) \oplus K_H(*) \right)^{\mathbb{Z}_2}.
}
\end{align}

Whenever $H \not \subset G_0$ we have that $H$ is a magnetic group
with $H_0 = G_0 \cap H$. Hence 
\begin{align}
\mathcal{K}_G(G/H) \cong \mathcal{K}_H(*) \ \ \ \mathrm{and} \ \ 
{K}_{G_0}(G/H) = K_{H_0}(*),
\end{align}
and the restriction homomorphism for the orbit type $G/H$ 
boils down to the restriction homomorphism for the magnetic group $H$:
\begin{align}
\xymatrix{
    \mathcal{K}_G(G/H) \ar[r] \ar[d]^\cong & K_{G_0}(G/H)^{\mathbb{Z}_2} \ar[d]^\cong \\
    \mathcal{K}_{H}(*) \ar[r] &  K_{H_0}(*) ^{\mathbb{Z}_2}.
}
\end{align}
Since the bottom horizontal map is an isomorphism rationally, then the upper horizontal map is also one. We conclude that the
the restriction homomorphism
\begin{align}
    \mathcal{K}_G(G/H) \to K_{G_0}(G/H)^{\mathbb{Z}_2}
\end{align}
induces an isomorphism rationally
\begin{align}
    \mathcal{K}_G(G/H) \otimes \mathbb{Q} \stackrel{\cong}{\to} K_{G_0}(G/H)^{\mathbb{Z}_2} \otimes \mathbb{Q}.
\end{align}

The previous argument can also be applied for $G$-spaces
of the form $G/H \times X$ where $X$ has no $G$-action, thus implying that the restriction homomorphism
\begin{align}
    \mathcal{K}_G(G/H \times X ) \to K_{G_0}(G/H \times X)^{\mathbb{Z}_2}
\end{align}
is an isomorphism rationally
\begin{align}
    \mathcal{K}_G(G/H \times X ) \otimes \mathbb{Q} \stackrel{\cong}{\to} K_{G_0}(G/H \times X)^{\mathbb{Z}_2} \otimes \mathbb{Q}.
\end{align}

The compatibility with open $G$-equivariant charts
on both sides of the restriction homomorphisms, the Mayer-Vietoris sequence, together with the five-lemma and the fact that the higher K-theory groups are defined with the usual suspension equality of eqn. \eqref{suspension}, implies that
the restriction homomorphism induces an isomorphism rationally.

$\square$

\subsection*{Twisted case}

The magnetic $G$-equivariant K-theory has an extension to the
case on which the local magnetic representations are projective. A simple way to understand this feature is the following.

Take $A \subset \mathbb{S}^1$ and consider an extension $\widetilde{G}$ of $G$ by $A$ fitting in the exact sequence
\begin{align}
   1 \to A \to \widetilde{G} \to G \to 1
\end{align}
where $G$ acts on $A$ by complex conjugation through the homomorphism $\phi$. That is $G \times A \to A$, $(g,a) \mapsto a^{-1}$ if $g \in G \backslash G_0$ and $(g,a) \mapsto a$ otherwise. Denote by $\widetilde{G}_0$ the group extension
over $G_0$ and note that we have the following diagram of exact sequences:
\begin{align}
    \xymatrix{
    A \ar[d] \ar@{=}[r] & A \ar[d] & \\
    \widetilde{G}_0 \ar[r] \ar[d] & \widetilde{G} \ar[r] \ar[d] & \mathbb{Z}_2 \ar@{=}[d] \\
    G_0 \ar[r] & G \ar[r] & \mathbb{Z}_2.
    }
\end{align}
The groups $\widetilde{G}$ and $G$ are magnetic, the former extending the latter in the middle vertical exact sequence, the left vertical exact arrow encodes the fact that $\widetilde{G}_0$
is a central $A$-extension of $G_0$, and the horizontal exact sequences encode the information of the magnetic groups $\widetilde{G}$ and $G$.

A {\it $\widetilde{G}$-twisted magnetic $G$-vector bundle}
over the $G$-space $X$ consists of a magnetic
$\widetilde{G}$-equivariant vector bundle $E \to X$ where the subgroup $A$
acts on the fibers of $E$ complex linearly by multiplication of scalars. Here $X$ is considered a $\widetilde{G}$-space by the induced action of the projection map $\widetilde{G} \to G$.

The {\it $\widetilde{G}$-twisted magnetic $G$-equivariant K-theory} of $X$, denoted as ${}^{\widetilde{G}}\mathcal{K}_G(X)$,
will be the subgroup of the magnetic $\widetilde{G}$-equivariant K-theory $\mathcal{K}_{\widetilde{G}}(X)$ generated
by $\widetilde{G}$-twisted magnetic $G$-vector bundles over $X$.

If we restrict the group $\widetilde{G}$ to $A$
\begin{align}
    \mathrm{res}^{\widetilde{G}}_A : \mathcal{K}_{\widetilde{G}}(X) \to K_A(X)
\end{align} and noting that $K_A(X) \cong R(A) \otimes K(X)$, we may take the irreducible representation $\nu$ of $A$ induced by the canonical inclusion $A \subset \mathbb{S}^1$.
Then the $\widetilde{G}$-twisted magnetic $G$-equivariant K-theory of $X$ can be understood as the following inverse image:
\begin{align}
    {}^{\widetilde{G}}\mathcal{K}_G(X) = (  \mathrm{res}^{\widetilde{G}}_A)^{-1} \left( \mathbb{Z}\langle \nu \rangle \otimes K(X) \right).
\end{align}

The higher K-theory groups ${}^{\widetilde{G}}\mathcal{K}^{-q}_G(X)$ are defined as in eqn. \eqref{suspension}, and therefore we obtain the inclusion:
\begin{align}
    {}^{\widetilde{G}}\mathcal{K}^{-q}_G(X) \subset \mathcal{K}^{-q}_{\widetilde{G}}(X).
\end{align}

The $\widetilde{G}$-twisted magnetic representations of $G$
define a subgroup of the magnetic representations of $\widetilde{G}$ and they 
split accordingly into real, complex and quaternionic type:
\begin{align}
    {}^{\widetilde{G}}\mathcal{R}(G) \subset \mathcal{R}(\widetilde{G}), \ \  \ \ \ {}^{\widetilde{G}}\mathcal{R}(G)_{\mathbb{F}} \subset \mathcal{R}(\widetilde{G})_{\mathbb{F}}, \ \ \ \mathrm{for} \ \ \ \mathbb{F} \in \{\mathbb{R}, \mathbb{C}, \mathbb{H} \}.
\end{align}

The coefficients of the twisted case ${}^{\widetilde{G}}\mathcal{K}^{-q}_G(*)$ split as described in eqn. \eqref{splitting suspension}:
\begin{align}
    {}^{\widetilde{G}}\mathcal{K}^{-q}_G(*)& \cong \\
   & \left( {}^{\widetilde{G}}\mathcal{R}(G)_{\mathbb{R}} \otimes KO^{-q}(*) \right) \oplus  \left({}^{\widetilde{G}}\mathcal{R}(G)_{\mathbb{C}}\otimes K^{-q}(*)  \right)  \oplus  \left( {}^{\widetilde{G}}\mathcal{R}(G)_{\mathbb{H}}\otimes KSp^{-q}(*) \right) \nonumber
\end{align}

The restriction to $\widetilde{G}_0$ gives us a natural homomorphism
\begin{align}
    {}^{\widetilde{G}}\mathcal{K}^{*}_G(X) \to {K}^{*}_{\widetilde{G}_0}(X)^{\mathbb{Z}_2}
\end{align}
inducing a rational isomorphism with its image. But note
that the image does not lie only in the $\widetilde{G}_0$-twisted complex K-theory 
${}^{\widetilde{G}_0}{K}^{*}_{G_0}(X)$ unless the group $A= \mathbb{Z}_2$. This follows from the fact that the conjugate representation of the canonical inclusion $ A \subset \mathbb{S}^1$ is only equal to the canonical one when the group $A$ is $\mathbb{Z}_2$.

Since the case $A=\mathbb{Z}_2$ is of independent interest, our main Theorem implies the following.

{\bf Corollary.} {\it Let $X$ be a compact $G$-space, $G \stackrel{\phi}{\to} {\mathbb{Z}_2}$ a magnetic group,
$\iota : G_0 \to G$ the inclusion of the kernel of $\phi$,
$\widetilde{G}$ a central extension of $G$ by the groups $\mathbb{Z}_2$ and $\widetilde{G}_0$ the kernel of the induced homomorphism from $\widetilde{G}$ to $\mathbb{Z}_2$.
Then the pullback of the restriction from $G$ to $G_0$ induces a homomorphism of twisted K-theories}
\begin{align}
    \iota^*: {}^{\widetilde{G}}\mathcal{K}_G^*(X) \to {}^{\widetilde{G}_0}K_{G_0}^*(X)^{\mathbb{Z}_2}
\end{align}
{\it and it becomes an isomorphism rationally}
\begin{align}
    \iota^* \otimes \mathbb{Q}: {}^{\widetilde{G}}\mathcal{K}_G^*(X) \otimes \mathbb{Q} \stackrel{\cong}{\to} {}^{\widetilde{G}_0}K_{G_0}^*(X)^{\mathbb{Z}_2}\otimes \mathbb{Q}.
\end{align}

\section{Applications}

The motivation for defining the magnetic equivariant K-theory groups comes from the realm of condensed matter physics.
The electronic properties of crystals, magnetic or not, can
be modeled with computers using Density Functional Theory (DFT) and the different programs that have been developed
for this task. The bundle of valence eigenvectors of the associated Hamiltonian, the Bloch bundle, in the case that the energy is gapped at the Fermi level, defines a magnetic
equivariant vector bundle. Adiabatic perturbations
of the Hamiltonian do not change the topological type of 
the Bloch bundle, and therefore the topological invariants
of the Bloch bundle, as an element in the magnetic equivariant
K-theory, become topological invariants of the Hamiltonian.

Several papers have been written on this regard (see \cite{Gomi} and references therein), but very
few explicit calculational tools have been developed in order
to extract the topological invariants from the Bloch bundle 
of a prescribed Hamiltonian. The main Theorem of this work
provides a calculational tool that may permit extract the non-torsion invariants of the Bloch bundle, in particular in the
case that the Hamiltonian models a magnetic material.

But before we see the applications in condensed matter physics, let us first start by relating the magnetic equivariant K-theory groups to other known K-theories.

When the magnetic group $G$ is $\mathbb{Z}_2$ and $\phi$ is the identity, the magnetic equivariant K-theory groups are exactly Atiyah's real K-theory groups \cite{Atiyah-KR}:
\begin{align}
    \mathcal{K}_{\mathbb{Z}_2}^* = K\mathbb{R}^*.
\end{align}
For spaces $X$ with trivial $\mathbb{Z}_2$ action we recover the K-theory of real vector bundles
\begin{align}
    \mathcal{K}_{\mathbb{Z}_2}^*(X) = KO^*(X).
\end{align}

When the magnetic group $G$ is again $\mathbb{Z}_2$ with $\phi$ the identity, and we choose its twisted version, namely $\widehat{G}= \mathbb{Z}_4$ and  $A \cong \mathbb{Z}_2$
acting by multiplication by $-1$, then we recover quaternionic K-theory \cite{Gomi}:
\begin{align}
    {}^{\mathbb{Z}_4}\mathcal{K}_{\mathbb{Z}_2}^* = K\mathbb{H}^*.
\end{align}
For spaces $X$ with trivial $\mathbb{Z}_2$ action we recover
the K-theory of symplectic vector bundles \cite{Dupont}:
\begin{align}
     {}^{\mathbb{Z}_4}\mathcal{K}_{\mathbb{Z}_2}^*(X) \cong KSp^*(X).
\end{align}

The well known relation between symplectic and real K-theory $KO^{*-4}=KSp^*$ also works for $\mathbb{Z}_2$-spaces and we have $\mathcal{K}_{\mathbb{Z}_2}^{*-4} \cong  {}^{\mathbb{Z}_4}\mathcal{K}_{\mathbb{Z}_2}^*$ \cite{Gomi}.
The restriction homomorphism $\iota^*$ lands in both cases in
complex K-theory $K^*$ and the induced $\mathbb{Z}_2$ action is simply the conjugation action on the complex vector bundles.
The generator of $K^{-2}$ is the virtual bundle $\mathcal{H}-1$
where $\mathcal{H}$ is the Hopf bundle over $S^2$ and $1$
is a trivial vector bundle of rank $1$. The conjugation map
sends $\mathcal{H}-1$ to $\overline{\mathcal{H}} - 1$
and in $K^{-2}$ one is the additive inverse of the other.

In $K^{-4}$ the generator is $\mathcal{P}-2$ where $\mathcal{P}$ is a rank $2$ complex vector bundle over $S^4$
whose clutching function  is given by the canonical diffeomorphism $\eta: S^3 \stackrel{\cong}{\to} SU(2)$.
Since the matrices of $SU(2)$ all commute with the matrix $i \sigma_y \mathbb{K}$ where $\sigma_y$ is the second Pauli matrix and $\mathbb{K}$ denotes complex conjugation, then it is clear that the conjugate clutching function $\overline{\eta}:= \mathbb{K} \eta \mathbb{K}$ is the same as the clutching function $(i \sigma_y) \eta (i \sigma_y)^{-1}$.
Therefore $\mathcal{P}$ and its conjugate $\overline{\mathcal{P}}$ are isomorphic. This implies that the conjugation action of $\mathbb{Z}_2$ is trivial on $K^{-4}$. By Bott periodicity we see that the complex K-theory groups
as $\mathbb{Z}_2$ modules are the following:
\begin{align}
K^{q} = \left\{
\begin{matrix}
    \mathbb{Z}  &\ \ \mathrm{trivial} \ &\ \mathbb{Z}_2-\mathrm{mod.}  \ & \ \mathrm{for} & \ q \equiv_4 0 \\ 
    \mathbb{Z}  &\ \ \mathrm{sign} \ &\ \mathbb{Z}_2-\mathrm{mod.} \ & \ \mathrm{for} & \ q \equiv_4 2 \\
    0 & \mathrm{for} & \ q \equiv_2 1. &&
\end{matrix}
\right.  \label{Z2 module K-theory}
\end{align}

We therefore have that the restriction maps
\begin{align}
\iota : \mathcal{K}_{\mathbb{Z}_2}^{*} \to (K^*)^{\mathbb{Z}_2} \ \ \ \ \iota :   {}^{\mathbb{Z}_4}\mathcal{K}_{\mathbb{Z}_2}^* \to (K^*)^{\mathbb{Z}_2} 
\end{align}
recover the well known rational isomorphisms:
\begin{align}
    \mathcal{K}_{\mathbb{Z}_2}^{*} \cong_{\mathbb{Q}}  {}^{\mathbb{Z}_4}\mathcal{K}_{\mathbb{Z}_2}^{*-4} \cong_{\mathbb{Q}}  (K^*)^{\mathbb{Z}_2} \cong_{\mathbb{Q}}  \left\{
    \begin{matrix}
        \mathbb{Q} & \mathrm{for} & * \equiv_4 0\\
        0 &  \mathrm{for} & * \not\equiv_4 0
    \end{matrix}
    \right.
\end{align}

Now that we have determined the $\mathbb{Z}_2$-module structure
of the complex K-theory groups in eqn. \eqref{Z2 module K-theory} we are now ready to show an application in condensed matter physics.

\subsection*{Topological Insulators Altermagnets}

Altermagnetism is a type of magnetic state in crystals
on which the magnetic structures are collinear and crystal-symmetry compensated resulting in zero net magnetization. But unlike ordinary collinear antiferromagnets,  the electronic bands in an altermagnets are not Kramer's degenerate \cite{altermagnet}. The symmetry that is present in some of
these 2D materials is $C_4\mathbb{T}$, a composition of a four fold rotation with time reversal symmetry.

Whenever the material is an insulator, various works have put forward the idea that there is a topological invariant that
separates trivial insulators from the topological ones \cite{C4T1,C4T2}.
It turns out that the magnetic equivariant K-theory predicts the existence of a bulk $\mathbb{Z}_2$-invariant on 2D systems
with the $C_4\mathbb{T}$ symmetry \cite{C4T-Hamiltonian}. This invariant can be extracted adding the spin $z$ operator $S_z$ on the system, and the reason for this procedure to work, lies in the rational isomorphism presented in the Corollary above.

The setup is as follows (for details we refer to  \cite{C4T-Hamiltonian}). The group is the one generated by $C_4\mathbb{T}$ and $S_z$ and they act on the 2D torus $T^2$. The operator $S_z$ squares to $1$ and acts trivially on $T^2$. The other operator acts on the space as follows:
\begin{align}
    C_4\mathbb{T} : T^2 \to T^2  \ \ \ (x,y) \mapsto (y,-x),
\end{align}
while on the fibers it incorporates the interaction of the spin with the lattice (spin orbit coupling). We have that $(C_4)^4=-1$ and $\mathbb{T}^2=-1$, therefore $(C_4\mathbb{T})^4=-1$.

The spin commutes with the rotation $C_4$ while it anticommutes with $\mathbb{T}$. Therefore $C_4 \mathbb{T}$ and $S_z$ anticommute on the fibers, while they clearly commute acting on $T^2$. To bring the notation we have employed in the previous section we have that the group of symmetries is $G=\mathbb{Z}_4 \times \mathbb{Z}_2 = \langle C_4 \mathbb{T} \rangle \times \langle S_z \rangle$, and the homomorphism $\phi$
is:
\begin{align}
\phi : \mathbb{Z}_4 \times \mathbb{Z}_2 \to \mathbb{Z}_2, \ \ \ \phi(a,b) = a  \ \mathrm{mod} \ 2.
\end{align}

Now, because of spin-orbit coupling, the group $G$ acts projectively on the fibers and therefore we need to take a central $A=\mathbb{Z}_2$ extension $\widetilde{G}$ of $G$, where
$A$ acts on the fibers by multiplication. We have the extension
\begin{align}
 \mathbb{Z}_2 \to \widetilde{G} \to G
\end{align}
where $\widetilde{G}$ is given by the following generators and relations
\begin{align}
\widetilde{G}= \langle \widetilde{a}, b | \widetilde{a}^8=b^2=1, bab=a^5 \rangle.
\end{align}

The equations above follow from the facts that $(C_4\mathbb{T})^4=-1$ and $S_z (C_4 \mathbb{T}) = -(C_4 \mathbb{T}) S_z = (C_4 \mathbb{T})^5S_z$. Here we have that
\begin{align}
    G_0= \langle a^2, b \rangle \cong \mathbb{Z}_2 \times \mathbb{Z}_2   \ \ \mathrm{and} \ \ \widetilde{G}_0= \langle \widetilde{a}^2, b \rangle \cong  \mathbb{Z}_4\times \mathbb{Z}_2.
\end{align}
Note that the restricted twisted group $\widetilde{G}_0$
becomes abelian. This is important for the calculation.

We want to determine the non-torsion invariants of the 
twisted equivariant magnetic group ${}^{\widetilde{G}}\mathcal{K}_G^*(T^2)$ and for this
we are going to calculate $ \left({}^{\widetilde{G}_0}K_{G_0}^*(T^2)\right)^{\mathbb{Z}_2}$. First note that
\begin{align}
    {}^{\widetilde{G}_0}K_{G_0}^*(T^2)\cong
 R(\mathbb{Z}_2)  \otimes {}^{\mathbb{Z}_4}K_{\mathbb{Z}_2}^*(T^2)  
\end{align}
since the spin $z$ operator $S_z$ commutes with $(C_4)^2=C_2$, and
$R(\mathbb{Z}_2)$ denotes the representation ring of $\mathbb{Z}_2$. The conjugation $\mathbb{Z}_2$-action generated by the operator $C_4 \mathbb{T}$ on 
the complex vector bundles splits into, 
conjugation composed the pullback of the four fold rotation in ${}^{\mathbb{Z}_4}K_{\mathbb{Z}_2}^0(T^2)$, while on $R(\mathbb{Z}_2)$ sends spin up to spin down and vice versa.

Therefore, as a $\mathbb{Z}_2$-module, we could take the isomorphism
\begin{align}
    {}^{\widetilde{G}_0}K_{G_0}^*(T^2) \cong \left( 
     {}^{\mathbb{Z}_4}K_{\mathbb{Z}_2}^*(T^2)  \right)^{\uparrow} \oplus 
      \left( 
     {}^{\mathbb{Z}_4}K_{\mathbb{Z}_2}^*(T^2)  \right)^{\downarrow}
\end{align}
where the arrows denote that we have split the bundle into spin up and spin down part. The $\mathbb{Z}_2$-action simply maps
spin up bundles to spin down, and vice versa. Hence, we could take the spin up bundles as the representatives of the $\mathbb{Z}_2$-invariants thus getting the isomorphism:
    \begin{align}
   \left( {}^{\widetilde{G}_0}K_{G_0}^*(T^2) \right)^{\mathbb{Z}_2} \cong \left[ \left( 
     {}^{\mathbb{Z}_4}K_{\mathbb{Z}_2}^*(T^2)  \right)^{\uparrow} \oplus 
      \left( 
     {}^{\mathbb{Z}_4}K_{\mathbb{Z}_2}^*(T^2)  \right)^{\downarrow} \right]^{\mathbb{Z}_2} \cong  
      \left( 
     {}^{\mathbb{Z}_4}K_{\mathbb{Z}_2}^*(T^2)  \right)^{\uparrow} 
\end{align}

Now, since $H^2(\mathbb{Z}_2, \mathbb{S}^1)=0$ we know
that $ {}^{\mathbb{Z}_4}K_{\mathbb{Z}_2}^*(T^2)$ is isomorphic
to the untwisted version $K_{\mathbb{Z}_2}^*(T^2)$ \cite{AdemRuan}. Therefore we conclude that the non-torsion invariants of the twisted equivariant magnetic K-theory that we were interested in are isomorphic to the non-torsion invariants of the complex K-theory $K_{\mathbb{Z}_2}^*(T^2)$:
\begin{align}
    {}^{\widetilde{G}}\mathcal{K}_G^*(T^2)  \cong_{\mathbb{Q}} 
    \left( 
     {}^{\mathbb{Z}_4}K_{\mathbb{Z}_2}^*(T^2)  \right)^{\uparrow}  \cong_{\mathbb{Q}} 
     K_{\mathbb{Z}_2}^*(T^2).
\end{align}

We can therefore conclude with the following result.

{\bf Result:} {\it Bundles over the torus $T^2$ with both the 
$C_4 \mathbb{T}$ symmetry and $S_z$ spin symmetry possess an integer invariant
coming from the 2D- cell (the bulk). This invariant can be extracted by determining the Chern number of the spin up bundle.}

The previous result was exploited in the work of the first two authors \cite{C4T-Hamiltonian} where it is furthermore shown
that there is a $\mathbb{Z}_2$ bulk invariant for systems that preserve the $C_4 \mathbb{T}$ symmetry on a 2D torus. Incorporating the spin $z$ it is shown that the value of this invariant can be  extracted as the parity of the Chern number of the spin up bundle defined above.

\section*{Conclusions}

We have shown that the torsion free part of the magnetic equivariant K-theory can be extracted from the conjugation invariant part of the restricted complex equivariant K-theory. We have used this result to show that the spin up Chern number of bundles with $C_4\mathbb{T}$ and spin $z$ symmetry over the 2D torus determine their bulk invariant. The parity of this invariant turns out to be the $\mathbb{Z}_2$ invariant that classifies topological insulators altermagnets with $C_4\mathbb{T}$ symmetry.

Further research is needed in order to determine the bulk invariants of magnetic symmetries that involve compositions of rotations, time reversal symmetry, and translations. These compositions require the definition of a more general kind of twistings that incorporate information of the base space. The associated twisted magnetic K-theories need to be defined and its properties determined. We leave this project for a future publication.

\section*{Acknowledgments}
H. Serrano acknowledges support from SECIHTI through the PhD scholarship No. 926934.
B. Uribe acknowledges the support of the Max Planck Institute for Mathematics in Bonn, Germany, and the continuous support of the Alexander Von Humboldt Foundation, Germany.
H. Serrano, B. Uribe and M. Xicot\'encatl acknowledge the financial support of SECIHTI grant CB-2017-2018-A1-S-30345.

The present paper is presented by B. Uribe to the Colombian Academy of Exact, Physical and Natural Sciences as part of the requirements to be appointed numerary member of the Academy. B. Uribe will succeed Academician Jaime Ignacio Lesmes Camacho (1939-2023) in the chair number 11. B. Uribe honors Jaime Lesmes as his professor, mentor and friend.

\section*{Authorship contribution statement}

All authors contributed equally on the results presented in the paper.

\section*{Conflict of interests}

The authors declare that they have no known competing financial interests or personal relationships that could have appeared to
influence the work reported in this paper.


    \bibliographystyle{alpha}
\bibliography{Bibliografia}
\end{document}